\documentclass[11pt]{amsart}
\usepackage[utf8]{inputenc}
\usepackage{amsmath,amssymb}
\usepackage{wrapfig}
\usepackage{url}
\usepackage{mathtools}
\usepackage{graphicx}
\usepackage{stmaryrd}
\usepackage{amsthm}
\usepackage{xcolor}
\usepackage[colorlinks=true,linkcolor=blue,citecolor=blue]{hyperref}
\usepackage[shortlabels]{enumitem}
\usepackage{comment}
\usepackage{relsize}
\usepackage{dsfont}
\usepackage{stmaryrd}
\usepackage{geometry}
\usepackage{setspace}
\usepackage{mathrsfs} 
\usepackage{orcidlink}
\usepackage{slashed}
\usepackage{yhmath}
\usepackage{amsfonts,amscd,bezier}
\usepackage{latexsym}
\usepackage[numbers, sort&compress]{natbib}
\makeatletter

\newtheorem*{theorem*}{Theorem}
\newtheorem{prop*}{Proposition}

\newtheorem{theorem}{Theorem}
\numberwithin{theorem}{section} 
\numberwithin{equation}{section}
\newtheorem{corollary}[theorem]{Corollary}

\newtheorem{prop}[theorem]{Proposition} 
\theoremstyle{definition}

\numberwithin{example}{section}

\newtheorem{lemma}[theorem]{Lemma}
\theoremstyle{definition}

\theoremstyle{definition}
\newtheorem{definition}[theorem]{Definition}

\allowdisplaybreaks

\newcommand{\R}{\mathbb{R}}
\newcommand{\C}{\mathbb{C}}
\newcommand{\N}{\mathbb{N}}
\newcommand{\Z}{\mathbb{Z}}
\newcommand{\T}{\mathbb{T}}
\title[Fourier analysis on a class of Schwartz spaces and applications] {Fourier analysis on a class of Schwartz spaces and applications}

\author[I. V. Coelho da Silva]{Isadora Vieira Coelho da Silva}
\address{Departamento de Matemática, Instituto de Ciências Matemáticas e de Computação, Universidade de São Paulo, São Carlos, SP, Brazil.\newline 
\hspace*{\parindent} Department of Mathematics: Analysis, Logic and Discrete Mathematics, Ghent University, Ghent, Belgium.}
\email{isadorasilva@usp.br}

\author[P. L. Dattori da Silva]{Paulo Leandro Dattori da Silva}
\address{Departamento de Matemática, Instituto de Ciências Matemáticas e de Computação, Universidade de São Paulo, São Carlos, SP, Brazil.}
\email{dattori@icmc.usp.br}

\author[M. Ruzhansky]{Michael Ruzhansky}
\address{Department of Mathematics: Analysis, Logic and Discrete Mathematics, Ghent University, Ghent, Belgium.\newline 
\hspace*{\parindent} School of Mathematical Sciences, Queen Mary University of London, United Kingdom.
}
\email{Michael.Ruzhansky@ugent.be}

\subjclass{Primary 42B05; Secondary 35A01, 35B65}

\keywords{Fourier series, homogeneous spaces,  Schwartz hypoellipticity, Schwartz solvability.}

\date{\today}

\begin{document}

\begin{abstract}
    The study of regularity and solvability of partial differential operators on non-compact manifolds requires the construction of adapted functional spaces. In this work we introduce $\mathcal{S}(\mathfrak{M}\times\mathbb{R}^n)$, the Schwartz space on the product of a homogeneous space $\mathfrak{M}$ and the Euclidean space $\mathbb{R}^n$. We characterize $\mathcal{S}(\mathfrak{M}\times\mathbb{R}^n)$ and its topological dual by employing discrete Fourier analysis in terms of eigenfunction expansions associated with a class of anisotropic elliptic operators on $\mathbb{R}^n$. We provide a characterization via partial Fourier coefficients and via full coefficients, that is, with respect to both components of $\mathfrak{M} \times \mathbb{R}^n$ simultaneously, which is achieved by studying the Fourier analysis on $\mathfrak{M}$. Our results provide new tools for the analysis of operators in this broader setting. As an application, we investigate the global existence and regularity of solutions for a class of partial differential equations on $\mathfrak{M}\times \mathbb{R}^n$.
\end{abstract}
\maketitle

\section{Introduction}
In order to study the solvability and regularity of solutions to partial differential operators on $\mathbb{R}^n$, it is natural to consider spaces of rapidly decaying functions, such as the Schwartz space. Recall that the Schwartz space $\mathcal{S}(\mathbb{R}^n)$ consists of all smooth functions $f \in C^\infty(\mathbb{R}^n)$ such that 
\[
    \sup_{x \in \mathbb{R}^n} |x^\beta \partial_x^\alpha f(x)| < \infty, \quad \forall \alpha, \beta \in \mathbb{N}_0^n.
\]
Recent works have also considered these classes, and more generally Gelfand–Shilov spaces, on the product manifold $\mathbb{T}^d \times \mathbb{R}^n$, where $\mathbb{T}^d$ denotes the $d$-dimensional torus; see, for example, \cite{AC2022, ACK2025, K2024, K2025}.
For $f \in C^\infty(\mathbb{T}^d \times \mathbb{R}^n)$ and for every  $\gamma \in \N_0^{d}$ and $\alpha, \beta \in\N_0^n$, define 
\[ 
    \|f\|_{\alpha, \beta, \gamma} \doteq \sup_{(t,x) \in \mathbb{T}^d \times \mathbb{R}^n} | x^\beta \partial_x^\alpha \partial_t^\gamma   f(t,x)|. 
\]
The space of Schwartz functions on $\T^d \times \R^n$ is given by
\[
\mathcal{S}(\mathbb{T}^d \times \mathbb{R}^n) \doteq \{ f \in C^\infty(\mathbb{T}^d \times \mathbb{R}^n); \ \|f\|_{\alpha, \beta, \gamma} < \infty, \ \forall \gamma \in \N_0^d \text{ and } \alpha, \beta \in\N_0^n\}.
\]
Equivalently,  $f \in  \mathcal{S}(\mathbb{T}^d \times \mathbb{R}^n)$ if and only if
\begin{equation}\label{eq sup L2}
    \sup_{t\in \T^d}\sum_{|\alpha+\beta| =s} \| x^\beta \partial_x^\alpha \partial_t^\gamma  f \|_{L^2(\R^n)} < \infty, 
\end{equation}
for every $s \in \N$ and $\gamma \in \N_0^d$. 
From \cite{K2025}, the topology on $\mathcal{S}(\mathbb{T}^d \times \mathbb{R}^n)$ is induced by the countable family of seminorms
\[
    p_N(f)\doteq \sum_{|\gamma|+|\alpha|+|\beta| \leq N} \|f\|_{\alpha, \beta, \gamma}, \quad N\in\N_0,
\]
that endows $\mathcal{S}(\mathbb{T}^d \times \mathbb{R}^n)$ with a Fréchet space structure.
Denote by $\mathcal{S}^\prime(\mathbb{T}^d \times \mathbb{R}^n)$ its dual space, i.e., the space of all linear functionals 
\[ 
    u: \mathcal{S}(\mathbb{T}^d \times \mathbb{R}^n)  \to \C
\]
such that there exist $N,C> 0$ satisfying
\begin{equation}\label{def S'}
    |\langle u,f\rangle | \leq C p_N(f), \quad \forall f \in \mathcal{S}(\mathbb{T}^d \times \mathbb{R}^n).
\end{equation}

An important tool in the study of regularity and solvability of partial differential operators in these spaces is the Fourier analysis. 
Such tool was
developed 
in \cite{K2025}
via the partial Fourier transform on $ \mathbb{T}^d \times \mathbb{R}^n$.  However, it is also useful to consider a discrete approach. In \cite{AC2022} and \cite{ACK2025}, the theory is constructed using a discrete basis of $L^2(\R^n)$ associated with a class of partial differential operators introduced in \cite{GPR2011}, in the context of Gelfand–Shilov spaces.
In the present work, we develop the corresponding theory for Schwartz space 
and its dual, associated with a broader class of partial differential operators presented in \cite{CGPR2019} and \cite{CRW2025}. More generally, we define the Schwartz spaces on $\mathfrak{M} \times \mathbb{R}^n$, where $\mathfrak{M}$ is a homogeneous space of the form
\[
    \mathfrak{M} \doteq H \backslash G = \{ Hg : g \in G \},
\]
for a compact Lie group $G$ and a closed subgroup $H\subset G$. In particular, when $H=\{e\}$, we recover $ \mathfrak{M} = G$, and thus our results include the case $\mathfrak{M}=\T^d$, and also of compact Lie groups.  Some properties of $\mathfrak{M}$ are presented later.

We first introduce the class of operators considered in this work.  Let $k,m \in 2\N$ and consider a positive partial differential operator $Q$ of the form
\begin{equation}\label{def Q}
    Q = Q(x,D) = \sum_{k|\alpha| + m|\beta| \leq km} c_{\alpha,\beta} x^\beta D_x^\alpha, \qquad D_x^\alpha = (-i)^{|\alpha|} \partial_x^\alpha,  
\end{equation}
with   complex coefficients $c_{\alpha,\beta} \in \mathbb{C}$. 
We also assume the anisotropic ellipticity condition (see \cite{CGPR2019, CRW2025}), i.e., 
\begin{equation}\label{anistropic eliptic}
    q_k(x,\xi) = \sum_{k|\alpha| + m|\beta| = km} c_{\alpha,\beta} x^\beta \xi^\alpha \neq 0, \quad \text{for} \quad (x,\xi) \in \R^n \times \Z^n \setminus \{(0,0)\},
\end{equation}
or equivalently, there exist $ C_1,  C_2 , R > 0 $ such that the full symbol
\[
    q(x, \xi) = \sum_{k|\alpha| + m|\beta| \leq km} c_{\alpha\beta} x^\beta \xi^\alpha
\]
satisfies 
\[
    C_1 \leq \frac{|q(x, \xi)|}{\Lambda_{k,m}(x, \xi)} \leq C_2, \quad |(x, \xi)| \geq R,
\]
where
\[
    \Lambda_{k,m}(x, \xi) = (1 + |x|^{2k} + |\xi|^{2m})^{1/2}, \quad (x, \xi) \in \mathbb{R}^{2n}.
\]

An important class of examples generalized by $Q$, studied in \cite{CDR2021}, is that of anharmonic oscillators given by
\begin{equation}\label{anharmonic}
    \mathcal{A}_{a,b} = (-\Delta)^a +|x|^{2b},
\end{equation}
for $a,b \in \N$. Although these operators do not always satisfy  the standard ellipticity assumptions required in, for instance, \cite{AC2022} and \cite{GPR2011}, they do satisfy the anisotropic ellipticity condition \eqref{anistropic eliptic}.

Under these conditions, according to \cite{CGPR2019}, there exists an orthonormal basis of  eigenfunctions $\{\varphi_j\}_{j\in\mathbb{N}} \subset \mathcal{S}(\mathbb{R}^n)$, with real eigenvalues, $\lambda_j >0$, such that $\lambda_j \to+ \infty$. Moreover, the eigenvalues have the following asymptotic behavior:
\begin{equation}\label{weyl asymptotic}
    \lambda_j \sim C j^{\frac{mk}{n(m+k)}}, \quad \text{as } j \to +\infty,
\end{equation}
for some positive constant $C$.
 In particular, any $u \in L^2(\mathbb{R}^n)$ (or $\mathcal{S}'(\mathbb{R}^n)$) can be expanded as
\[
    u = \sum_{j =1}^{+\infty} u_j \varphi_j,
\]
where  the coefficients $u_j\in \C$ are
defined by $u_j = \langle u, \varphi_j  \rangle_{L^2(\mathbb{R}^n)}$ (or simply $u_j = \langle u, \varphi_j \rangle$) for all $j\in\N$.

In \cite{CRW2025}, one also considers the Hilbert subspaces $M_{\sigma,s}(\mathbb{R}^n)$ of $\mathcal{S}'(\mathbb{R}^n)$, for $\sigma > 0$ and $s \in \mathbb{R}$, defined by the norm
\[
    \|u\|_{M_{\sigma,s}} = \left( \iint_{\mathbb{R}^{2n}} |V_\varphi u(x, \xi)|^2 \, \theta_\sigma(x, \xi)^{2s} \, dx \, d\xi \right)^{\frac{1}{2}},
\]
where 
\[
    \theta_{\sigma}(x, \xi) = 1 + |x| + |\xi|^{\frac{1}{\sigma}}, \quad x, \xi \in \mathbb{R}^n,
\]
and $V_\varphi$ denotes the short-time Fourier transform of a tempered distribution $ u \in \mathcal{S}'(\mathbb{R}^n) $ defined by
\[
    V_{\varphi} u(x, \xi) = (2\pi)^{-\frac{n}{2}} \langle u, M_\xi T_x \varphi \rangle , \quad x, \xi \in \mathbb{R}^n, \quad \varphi \in \mathcal{S}(\mathbb{R}^n) \setminus \{0\},
\]
with 
\[ 
    T_x f(y) = f(y - x) \quad \text{and} \quad M_\xi f(y) = e^{i\langle y, \xi \rangle} f(y), \  x, y, \xi \in \mathbb{R}^n,
\]
for $ f $ a function or distribution defined on $ \mathbb{R}^n $.
For any $ \sigma > 0 $, $ M_{\sigma,0}(\mathbb{R}^n) = L^2(\mathbb{R}^n) $, and $ M_{\sigma,s_1}(\mathbb{R}^n) \subseteq M_{\sigma,s_2}(\mathbb{R}^n) $ whenever $ s_1 \geq s_2 $. Moreover,
\begin{equation}\label{S S'}
    \mathcal{S}(\mathbb{R}^n) = \bigcap_{s \in \mathbb{R}} M_{\sigma,s}(\mathbb{R}^n) \quad \text{and} \quad \mathcal{S}'(\mathbb{R}^n) = \bigcup_{s \in \mathbb{R}} M_{\sigma,s}(\mathbb{R}^n). 
\end{equation}

By Proposition 3.1 of \cite{CRW2025}, for $ \sigma = \frac{k}{m} $  we have
\begin{equation}\label{Msigma}
    \|u\|_{M_{\frac{k}{m},s}(\mathbb{R}^n)}^2 \sim \sum_{j=1}^{+\infty} \lambda_j^{\frac{2s}{k}} |u_j|^2 = \| Q^{\frac{s}{k}} u \|^2_{L^2(\R^n)}, \quad \forall s>0,
\end{equation}
where the operator $ Q^s $, for $s\in \R$, is  defined by
\begin{equation*}\label{def Qs}
    Q^s u = \sum_{j=1}^{+\infty} \lambda_j^s u_j \varphi_j, \quad \forall  u \in \mathcal{S}'(\mathbb{R}^n).
\end{equation*}
Also, from Corollary 3.4 of \cite{CRW2025},  the Schwartz space $\mathcal{S}(\mathbb{R}^n)$ is characterized in terms of the coefficients with respect to the basis $\{\varphi_j\}$, as
\begin{equation}\label{Schwartz jN}
    \begin{aligned}
        u\in \mathcal{S}(\mathbb{R}^n) &\iff |u_j| = O(j^{-N}), \quad j\to +\infty, \ \forall N \in \mathbb{N},\\
        &\iff |u_j| = O(\lambda_j^{-N}), \quad j\to +\infty, \ \forall N \in \mathbb{N}.
    \end{aligned}
\end{equation}

In this work, we extend this characterization to the spaces $\mathcal{S}(\mathfrak{M} \times \mathbb{R}^n)$ and $\mathcal{S}'(\mathfrak{M} \times \mathbb{R}^n)$ in terms of partial Fourier coefficients associated with $\{\varphi_j\}_{j \in \mathbb{N}}$. Moreover, we provide a characterization via full coefficients, that is, with respect to both components of $\mathfrak{M} \times \mathbb{R}^n$ simultaneously. This is achieved using  Fourier analysis on $\mathfrak{M}$, which we describe in the next section. In the last section, as an application, we study the Schwartz regularity and solvability of a class of operators on  $\mathfrak{M}\times \R^n$, given by
\[ 
    \mathsf{P} = T + c Q + q, 
\]
 where $c, \ q \in \C$, $c\neq 0$, and $T$ is a left-invariant symmetric operator on $\mathfrak{M}$.

\section{Homogeneous spaces}

Let $G$ be a compact Lie group of dimension $d$, and let $H \subset G$ be a closed subgroup. In \cite{RT2010}, the compact homogeneous space $\mathfrak{M}$ is defined as the set of classes $Hg = \{h g : h \in H \}$ for every $g \in G$, that is,
\[
    \mathfrak{M} \doteq H \backslash G = \{ Hg : g \in G \}.
\]
We endow $\mathfrak{M}$  with the topology co-induced by $g \mapsto Hg$,  equipped with the unique $ C^\infty$-manifold structure such that the mapping 
\begin{align*}
    G\times H \backslash G &\to H \backslash G\\
    (g_1, H g_2) &\mapsto H g_2 g_1
\end{align*}
is smooth and such that there exist a neighborhood $ U \subset H \backslash G $  of $ H e \in H\backslash G $  and a mapping $ \psi \in C^\infty(U, G) $ satisfying $H\psi(Hg) = Hg$ for $Hg \in U$ (Definition 13.2.2 in \cite{RT2010} and Theorem 2.1.7 in \cite{C2014}).

The properties of function spaces on $\mathfrak{M}$ are closely related to those on $G$. Following \cite{M2020, RT2010}, we now recall some useful properties of $G$.

Denote  by $\mu_G$ the normalised Haar measure on $G$, i.e., the unique measure $\mu_G$  that satisfies 
\[
    \int_G 1 \, d\mu_G = \mu_G(G) = 1 \quad \text{and} \quad  \int_G f(g) \, d\mu_G(g) = \int_G f(yg) \, d\mu_G(g), \ \ \forall \ y \in G.
\]
Let $\mbox{Rep}(G)$ be the set of continuous irreducible unitary representations of $G$. Since $G$ is compact, every $\eta\in \mbox{Rep}(G)$ is finite-dimensional, and  it can be viewed as a matrix-valued function $\eta: G \to \mathbb{C}^{d_\eta \times d_\eta}$, where $d_{\eta}=\mathrm{dim} \ \eta$. We denote by $\widehat{G} $ the quotient of $\mbox{Rep}(G)$ by the equivalence relation given by: $\eta \sim \psi $ if there exists a  unitary
matrix $A \in \mathbb{C}^{d_\eta \times d_\eta} $ such that $A\eta(x) = \psi(x),$ for all $x\in G$.
By the Peter-Weyl Theorem (Theorem 7.5.14 of \cite{RT2010}), 
\[
    \mathcal{B}\doteq \{\sqrt{d_{\eta}}\eta_{rs}\ | \ 1 \leq r,s \leq d_{\eta}, \eta=(\eta_{rs})_{r,s=1}^{d_{\eta}}\in[\eta]\in \widehat{G}\},
\]
is an orthonormal basis for $L^{2}(\mu_{G})$, taking only one representative from each equivalence class.

Let $\{X_1, \dots, X_d\}$ be a basis of the Lie algebra $\mathfrak{g}$ of $G$.  
For a multi-index $\gamma = (\gamma_1, \dots, \gamma_d) \in \mathbb{N}_0^d$, consider the left-invariant differential operator of order $|\gamma| = \gamma_1 + \cdots + \gamma_d$ on $G$ by
\begin{equation}\label{def partialg}
    \partial_g^\gamma = Y_1 Y_2 \cdots Y_{|\gamma|},
\end{equation}
where each $Y_j \in \{X_1, \dots, X_d\}$ and for each $k = 1, \dots, d$, $ \displaystyle\sum_{j : Y_j = X_k  } 1 =\gamma_k$.
We denote by $\mathcal{L}_{G}$ the Laplace-Beltrami operator on $G$, that can be written as 
\[
    \mathcal{L}_G = - \sum_{i=1}^d X_i^2.
\]
For each $[\eta] \in \widehat{G}$ there exists $\nu_{[\eta]} \geq 0$ (independent of the choice of the representative) such that
\[
    -\mathcal{L}_{G} \eta(x)_{rs}=\nu_{[\eta]} \eta(x)_{rs}\ \text{for all} \ 1 \leq r,s \leq d_{\eta}.
\]
We use the notation $\langle\eta\rangle\doteq\left(1+\nu_{[\eta]}\right)^{1/2}$ for  the eigenvalues of $(I-\mathcal{L}_{G})^{1/2}$. For all non-trivial $[\eta] \in \widehat{G}$, there exists a constant $C>0$ such that
\begin{equation}\label{<eta>^2}
    \nu_{[\eta]} \leq \langle \eta \rangle^2 \leq C \nu_{[\eta]}.
\end{equation}
Moreover, for all $[\eta] \in \widehat{G}$,
\begin{equation}\label{eq deta}
    d_\eta \leq C \langle \eta \rangle^{\frac{d}{2}}
\end{equation}
and, for every $\gamma \in \N_0^d$ there exist $C,M>0$ such that
\begin{equation}\label{ex 1.39 [W]}
    |\partial_g^ \gamma \eta(g)_{rs}|\leq C\langle \eta \rangle^{M}, \quad \forall\  1\leq r,s\leq d_\eta,
\end{equation}
for all $g\in G$ (see Example 1.39 of \cite{M2020}). Also,  
\begin{equation}\label{sum deta < infty}
    \sum_{[\eta] \in \widehat{G}} d_\eta^2 \langle \eta \rangle^{-2t} < + \infty \iff t > \frac{d}{2}.
\end{equation}

For $f \in L^1(G)$, and  $[\eta] \in \widehat{G}$, $\eta = (\eta_{rs})_{r,s=1}^{d_\eta}$, the $\eta$-Fourier coefficient of $f$ is given by
\[
    \hat{f}(\eta) = \int_G f(x)\eta(x)^* d\mu_G \in \mathbb{C}^{d_\eta \times d_\eta},
\]
or, more explicitly,
\[
    \hat{f}(\eta)_{rs} = \int_G f(x)\overline{\eta(x)_{sr}} d\mu_G = \langle f, \eta_{sr} \rangle_{L^2(G)}.
\]
By the Peter-Weyl Theorem, we can write  $f \in L^2(G)$ by the Fourier series
\begin{align*}
    f(x) &= \sum_{[\eta] \in \widehat{G}} d_\eta \sum_{r,s=1}^{d_\eta} \langle f, \eta_{rs} \rangle_{L^2(G)} \eta(x)_{rs} = \sum_{[\eta] \in \widehat{G}} d_\eta \sum_{r,s=1}^{d_\eta} \hat{f}(\eta)_{sr} \eta(x)_{rs}. 
\end{align*}
The space $C^\infty(G)$ can be characterized through the Fourier coefficients as follows.

\begin{theorem}\label{thm C^infty(G)}
    Let $G$ be a compact Lie group. Then, $f\in C^\infty(G)$ if and only if for each $N>0$, there exists $C_N> 0$ such that 
    \[
        |\widehat{f}(\eta)_{rs}| \leq C_N\langle\eta\rangle^{-N},
    \] 
    for all $[\eta]\in\widehat{G}$, $1\leq r,s\leq d_\eta$.
\end{theorem}

We consider the space of distributions $\mathcal{D}'(G)$ as the space of all continuous linear functionals on $C^\infty(G)$, with the usual convergence: for $u_j, u \in \mathcal{D}'(G)$,  we say that $u_j \to u$ in $\mathcal{D}'(G)$ as $j \to +\infty$ if $u_j(\varphi) \to u(\varphi)$ in $\mathbb{C}$ as $j \to + \infty$, for all $\varphi \in C^\infty(G)$. We denote $\langle u, \varphi \rangle_G := u(\varphi)$, for all  $\varphi \in C^\infty(G)$.

For every $[\eta] \in \widehat{G} $, $\eta = (\eta_{rs})_{r,s=1}^{d_\eta}$, the $\eta$-Fourier coefficient of $u \in \mathcal{D}'(G)$ is given by
\[
    \widehat{u}(\eta) = \langle u, \eta^* \rangle_G \in \mathbb{C}^{d_\eta \times d_\eta},
\]
or, more explicitly,
\[
    \widehat{u}(\eta)_{rs} = \langle u, \overline{\eta_{sr}} \rangle_G.
\]    
We can also characterize the space $\mathcal{D}'(G)$ as follows.

\begin{theorem}\label{thm D'(G)}
    Let $G$ be a compact Lie group. Then, $u\in \mathcal{D}'(G)$ if and only if there exist $C,N>0$ such that 
    \[
        |\widehat{u}(\eta)_{rs}| \leq C \langle\eta\rangle^{N},
    \]    
    for all $[\eta]\in\widehat{G}$, $1\leq r,s\leq d_\eta$.
\end{theorem}

We are now in a position to present the properties of the space $\mathfrak{M}$, which are described in \cite{C2014, RT2010}.
For each $\eta \in \operatorname{Rep}(G)$, consider the space  $H_\eta$ of $\eta$-invariant vectors relative to $H$
given by
\[
    H_\eta \doteq \{ v \in \mathbb{C}^{d_\eta} : \eta(h)v = v \text{ for all } h \in H \}.
\]
Let $\widehat{G}_0 \doteq \{[\eta] \in \widehat{G}: \ H_\eta \neq \{0\}\}$ and $h_\eta \doteq \dim H_\eta$. Fix an orthonormal basis of $\mathbb{C}^{d_\eta}$ such that its first $h_\eta$ vectors form a basis of $H_\eta$.

For $f \in C^\infty(\mathfrak{M})$, we denote by $\widetilde{f} \in C^\infty(G)$ its lifting $\widetilde{f}(g) \doteq f(Hg)$. We endow $\mathfrak{M}$ with the measure $\mu_\mathfrak{M}$ given by the projection of $\mu_G$, satisfying
\[
    \int_\mathfrak{M} f  \, d\mu_\mathfrak{M} =  \int_G \tilde{f}  \, d\mu_G,
\]
and $ \| f\|_{L^2(\mu_\mathfrak{M})} = \| \tilde{f} \|_{L^2(\mu_G)}  $, for all $f\in C^\infty(\mathfrak{M})$. Also, if  $X\in \mathfrak{g}$, we define $X f \doteq X \tilde{f}$, for all $f\in C^\infty(\mathfrak{M})$.

Consider a continuous linear operator  $P: C^\infty(G) \to C^\infty(G) $, whose symbol is given by
\[
    \sigma_P(g,\eta) \doteq \eta^* (P\eta)(g) \in \C^{d_\eta\times d_\eta}.
\]
If $ \sigma_P(hg,\eta) = \sigma_P(g,\eta) , $
for all $g\in G$, $h\in H$ and $[\eta]\in \widehat{G}$, then $P$ maps $C^\infty(\mathfrak{M}) $ into itself. 
For instance, this holds
when $P$ is a left invariant operator, since the symbol does not depend on $g\in G$.

From Theorem 2.6.11 of \cite{C2014}, we can characterize the elements of $C^\infty(\mathfrak{M})$ in terms of their Fourier coefficients in $C^\infty(G) $. Then, for $f \in C^\infty(G)$ (or $L^2(G)$), we have $f \in C^\infty(\mathfrak{M})$ (or $L^2(\mathfrak{M})$) if and only if
\[
    \hat{f}(\eta)_{sr} = 0 \text{ for } [\eta] \notin \widehat{G}_0 \text{ or } r > h_\eta.
\]
Moreover, for any distribution $f\in \mathcal{D}'(\mathfrak{M}) \subset \mathcal{D}'(G)$, $\hat{f}(\eta) =0 $ for $[\eta]\notin \widehat{G}_0$,  $\hat{f}(\eta)_{sr} =0 $ for $r>h_\eta$, and we adopt the convention $ \eta(g)_{rs} = 0$ for all $r>h_\eta$. Then we write the Fourier series of $f$ by
\[
    f(g) = \sum_{[\eta] \in \widehat{G}_0} d_\eta \sum_{s=1}^{d_\eta} \sum_{r=1}^{h_\eta} \widehat{f}(\eta)_{sr} \, \eta(g)_{rs}.
\]

\section{\texorpdfstring{Fourier Analysis on $\mathfrak{M} \times \mathbb{R}^n$}{Fourier Analysis on MxRn}}

Let $f$ be a function in $ C^\infty(\mathfrak{M} \times \mathbb{R}^n)$.The lifting  $\widetilde{f} \in C^\infty(G \times \mathbb{R}^n)$ is given by
\[
    \widetilde{f}(g, x) \doteq f(Hg, x), \quad \forall (g, x) \in G \times \mathbb{R}^n.
\]
For multi-indices $\alpha, \beta \in \mathbb{N}_0^n$ and $\gamma \in \mathbb{N}_0^d$, define the seminorm
\[
    \|f\|_{\alpha, \beta, \gamma} \doteq \sup_{(g, x) \in G \times \mathbb{R}^n} \left| x^\beta \, \partial_x^\alpha \, \partial_g^\gamma \widetilde{f}(g, x) \right|,
\]
where $\partial_x^\alpha = \partial_{x_1}^{\alpha_1} \cdots \partial_{ x_n}^{\alpha_n}$, $x^\beta = x_1^{\beta_1} \cdots x_n^{\beta_n}$ and $ \partial_g^\gamma$ is described in \eqref{def partialg}. 
In what follows, we identify functions on $\mathfrak{M}$ with their liftings to $G$.

\begin{definition}
    The {Schwartz space} on $\mathfrak{M} \times \mathbb{R}^n$ is defined by
    \[
        \mathcal{S}(\mathfrak{M}  \times \mathbb{R}^n) \doteq \left\{ f \in C^\infty(\mathfrak{M} \times \mathbb{R}^n) \;:\; \|f\|_{\alpha, \beta, \gamma} < \infty \;\; \forall \alpha, \beta \in \mathbb{N}_0^n,\; \forall \gamma \in \mathbb{N}_0^d \right\}. 
    \]
\end{definition}

As in the case of the torus, we consider on $\mathcal{S}(\mathfrak{M} \times \mathbb{R}^n)$ the topology induced by the countable family of seminorms 
\[
    p_N(f)\doteq \sum_{|\gamma|+|\alpha|+|\beta| \leq N} \|f\|_{\alpha, \beta, \gamma}, \quad N\in\N_0.
\]    
The dual space  $\mathcal{S}^\prime(\mathfrak{M} \times \mathbb{R}^n)$ is the space of all linear functionals
\[
    u: \mathcal{S}(\mathfrak{M} \times \mathbb{R}^n) \to \mathbb{C}
\]
for which there exist $N, C > 0$ such that
\begin{equation}\label{def S' M}
    |\langle u,f\rangle| \leq C p_N(f), \quad \forall f \in \mathcal{S}(\mathfrak{M} \times \mathbb{R}^n).
\end{equation}

The space $\mathcal{S}(\mathfrak{M} \times \mathbb{R}^n)$ can be equivalently defined as follows.

\begin{prop}
    Let ${f} \in C^\infty(\mathfrak{M} \times \mathbb{R}^n)$. Then $f \in \mathcal{S}(\mathfrak{M} \times \mathbb{R}^n)$ if and only if
    \begin{equation}\label{eq sup L2 M}
        \sup_{g \in G} \sum_{|\alpha+\beta| = s} \| x^\beta \partial_x^\alpha \partial_g^\gamma {f}(g,\cdot) \|_{L^2(\mathbb{R}^n)} < +\infty,
    \end{equation}
    for every $s \in \mathbb{N}$ and $\gamma \in \mathbb{N}_0^d$.
\end{prop}
\begin{proof}
    Assume $f \in \mathcal{S}(\mathfrak{M} \times \mathbb{R}^n)$. Then, for all $\alpha,\beta \in \N_0^n$ and $\gamma \in \N_0^d$, there exists $C_{\alpha,\beta,\gamma}>0$ such that
    \[
        \|f\|_{\alpha,\beta,\gamma} = \sup_{(g,x) \in G \times \mathbb{R}^n} |x^\beta \partial_x^\alpha \partial_g^\gamma f(g,x)| < C_{\alpha,\beta,\gamma}.
    \]
    For $N\in \N$  and $1\le i\le n$, let  $\beta^i =(\beta_1^i, \dots, \beta_n^i) $ be given by $\beta^i_i = N$ and $\beta^i_j =0$ if $j\neq i$. For  $x=(x_1,\dots,x_n)\neq 0 $, let $i_0$ be such that  $|x_{i_0}| = \displaystyle\max_{1\leq i\leq n}\{|x_i | \} \neq 0$.
    Then
    \[ 
        |\partial_x^\alpha \partial_g^\gamma f(g,x)| \leq \dfrac{1}{|x_{i_0}|^N} C_{\alpha,\beta^{i_0},\gamma} \leq \max_{1\leq i\leq n}C_{\alpha,\beta^{i},\gamma} \dfrac{(\sqrt{n})^N}{|x|^N} = C_{\alpha, N,\gamma} |x|^{-N},
    \]
    where $C_{\alpha, N,\gamma} \doteq \displaystyle \max_{1\leq i\leq n}C_{\alpha,\beta^{i},\gamma} {(\sqrt{n})^N}$. 
    Fix $s \in \mathbb{N}$ and $\gamma \in \mathbb{N}_0^d $. Let $\alpha,\beta \in \N_0^n$ be such that $|\alpha+\beta| = s$. 
    We have  
    \begin{align*}
        \| x^\beta \partial_x^\alpha \partial_g^\gamma f(g,\cdot) \|_{L^2(\mathbb{R}^n)}^2 &= \int_{\R^n} |x^\beta \partial_x^\alpha \partial_g^\gamma f(g,x)|^2dx\\
        &= \int_{|x|\leq 1} |x^\beta \partial_x^\alpha \partial_g^\gamma f(g,x)|^2dx + \int_{|x|\geq 1} |x^\beta \partial_x^\alpha \partial_g^\gamma f(g,x)|^2dx\\
        &\leq C'_{\alpha,\beta,\gamma}  + \int_{|x|\geq 1} |x^{\beta}|^2 | \partial_x^\alpha \partial_g^\gamma f(g,x)|^2 dx\\
        &\leq C'_{\alpha,\beta,\gamma}  + C_{\alpha, N,\gamma}^2 \int_{|x|\geq 1} |x|^{2|\beta| - 2 N} dx 
    \end{align*}
   Taking $N$ sufficiently large, the last integral is finite.  Hence,
   \[
        \sup_{g \in G} \| x^\beta \partial_x^\alpha \partial_g^\gamma f(g,\cdot) \|_{L^2(\mathbb{R}^n)} < +\infty, 
   \]
   which implies that \eqref{eq sup L2 M} holds.

   On the other hand, assume \eqref{eq sup L2 M} holds for every $s \in \mathbb{N}$ and $\gamma \in \mathbb{N}_0^d$. Fix $\alpha_0, \beta_0 \in \N_0^n$ and $\gamma_0 \in \N_0^d$. We need to show that
    \[
        \sup_{(g,x) \in G \times \mathbb{R}^n} |x^{\beta_0} \partial_x^{\alpha_0} \partial_g^{\gamma_0} f(g,x)| < \infty.
    \]    
    For a fixed $g \in G$, define $h_g(x) \doteq x^{\beta_0} \partial_x^{\alpha_0} \partial_g^{\gamma_0} f(g,x)$. By Sobolev inequalities (see, for example,  Corollary 9.13 in \cite{B2010}), there exist an integer $m > n/2$ and a constant $C > 0$ (independent of $g$) such that
    \begin{align*}
        \sup_{x \in \mathbb{R}^n} |h_g(x)| &\leq C \sum_{|\delta| \le m} \| \partial_x^\delta h_g \|_{L^2(\mathbb{R}^n)}\\
        &\leq C \sum_{|\delta| \le m} \sum_{\gamma \leq \min\{\beta_0,\delta \} } C_\gamma \| x^{\beta_0 - \gamma} \partial_x^{\alpha_0 + \delta - \gamma} \partial_g^{\gamma_0} f(g,\cdot) \|_{L^2(\mathbb{R}^n)}.
    \end{align*}
    By hypothesis, the last sum is finite. Therefore, taking the supremum over $g \in G$, we conclude the proof. 
\end{proof}

Consider a positive partial differential operator $Q$ given by \eqref{def Q}, satisfying the anisotropic ellipticity condition \eqref{anistropic eliptic}, and let $\{\varphi_j\}_{j \in \mathbb{N}} \subset \mathcal{S}(\mathbb{R}^n)$ be the basis generated by $Q$, with eigenvalues $\{\lambda_j\}_{j \in \mathbb{N}}$ satisfying \eqref{weyl asymptotic}.

We first characterize $\mathcal{S}(\mathfrak{M} \times \mathbb{R}^n)$ in terms of the partial Fourier coefficients in the second variable. This result is an adaptation of results for Gelfand-Shilov spaces on $\T^1 \times \R^n$ (see Theorem 2.4 in \cite{AC2022}). 

\begin{theorem}\label{thm partial schwartz M}
    Let $u\in \mathcal{S}'(\mathfrak{M} \times \mathbb{R}^n)$. Then,  $u \in\mathcal{S}(\mathfrak{M} \times \mathbb{R}^n)$  if and only if it can be represented as
    \[
        u(g,x) = \sum_{j=1}^{\infty} u_j(g) \varphi_j(x),
    \]
    with $u_j \in C^\infty(\mathfrak{M})$ such that  for every $\gamma \in \mathbb{N}_0^d$ and  $N \in \mathbb{N}$, there exists a constant $C_{\gamma,N} > 0$ satisfying
    \begin{equation}\label{eq partial coef uj}
        \sup_{g \in G} |\partial_g^\gamma {u}_j(g)| \leq C_{\gamma,N} \, j^{-N}, \quad \forall j \in \mathbb{N}.
    \end{equation}
\end{theorem}

In order to prove this theorem, we first prove the following result. 

\begin{prop}\label{prop 2.3 [AC]22 M}
    Let $Q$ be the operator given by \eqref{def Q} and satisfying \eqref{anistropic eliptic}. Then, $u \in \mathcal{S}(\mathfrak{M} \times \mathbb{R}^n)$ if and only if for every $M \in \mathbb{N}$ and $\gamma \in \mathbb{N}_0^d$
    \begin{equation}\label{eq QM M}
        \sup_{g \in G} \|Q^M \partial_g^\gamma {u}(g,\cdot)\|_{L^2(\mathbb{R}^n)} < \infty.
    \end{equation}
\end{prop}
\begin{proof}
    Assume that $u \in \mathcal{S}(\mathfrak{M} \times \mathbb{R}^n)$, and let $M\in \N$. Since $Q^M$ can be written in the form
    \[
        Q^M = \sum_{k|\alpha|+m|\beta| \leq kmM} c'_{\alpha,\beta} x^\beta D_x^\alpha
    \]
    for some $c'_{\alpha,\beta} \in \mathbb{C}$, it follows from \eqref{eq sup L2 M} that
    \begin{align*}
        \|Q^M \partial_g^\gamma {u}(g,\cdot)\|_{L^2(\mathbb{R}^n)} &\leq \sum_{s=0}^{kmM} \sum_{k|\alpha|+m|\beta| =s} |c'_{\alpha,\beta}| \|x^\beta \partial_x^\alpha \partial_g^\gamma {u}(g,\cdot)\|_{L^2(\mathbb{R}^n)}\\
        &\leq C^M \sum_{s=0}^{kmM} \sum_{|\alpha|+|\beta| =s} \|x^\beta \partial_x^\alpha \partial_g^\gamma {u}(g,\cdot)\|_{L^2(\mathbb{R}^n)} < +\infty,
    \end{align*}
    for all $g \in G$, where we are using $\displaystyle\sup_{k|\alpha|+m|\beta|\leq kmM} |c'_{\alpha,\beta}| < C^M$ for some $C>0$.
    
    On the other hand, suppose that \eqref{eq QM M} holds for every $M \in \mathbb{N}$ and $\gamma \in \mathbb{N}_0^d$. In particular, for each fixed $g \in G$ and $\gamma \in \mathbb{N}_0^d$, from \eqref{Msigma} we conclude that 
    \[
        \| \partial_g^\gamma {u}(g,\cdot)\|_{M_{\sigma,s}} < +\infty, \quad \forall s = kM,\ M \in \mathbb{N}, \ \sigma=\dfrac{k}{m}.
    \]
    Therefore, since $M_{\sigma,s_1}(\mathbb{R}^n) \subseteq M_{\sigma,s_2}(\mathbb{R}^n)$ whenever $s_1 \geq s_2$, by \eqref{S S'} we have $\partial_g^\gamma {u}(g,\cdot) \in \mathcal{S}(\mathbb{R}^n)$ for each $g \in G$.
    
    Denote
    \[
        |\partial_g^\gamma {u}(g,\cdot)|_s \doteq \sum_{|\alpha+\beta| = s} \|x^\beta \partial_x^\alpha \partial_g^\gamma {u}(g,\cdot)\|_{L^2(\mathbb{R}^n)}.
    \]
    By Proposition 4.1 from \cite{GPR2011}, there exists $C > 0$ such that for every $s \in \mathbb{N}$ with $s = pm + r$, $p \in \mathbb{N}$, $0 < r < m$ and for all $\epsilon > 0$,
    \[
        |v|_s \leq \epsilon |v|_{(p+1)m} + C \epsilon^{-\frac{r}{m-r}} |v|_{pm} + C^s (s!)^{\frac{1}{2}} \|v\|_{L^2(\mathbb{R}^n)},
    \]
    for every $v \in \mathcal{S}(\mathbb{R}^n)$. Therefore, by \eqref{eq sup L2 M}, it remains to prove that for all $p,m \in \mathbb{N}_0$ there exists $C_{pm,\gamma} > 0$ such that
    \[
        \sup_{g \in G} |\partial_g^\gamma {u}(g,\cdot)|_{pm} \leq C_{pm,\gamma}.
    \]
    Let $\delta > 0$, $\gamma \in \mathbb{N}_0^d$ and denote
    \[
        \sigma_p(\partial_g^\gamma {u}, \delta) = (pm)!^{-\frac{km}{k+m}} \delta^{-p} |\partial_g^\gamma {u}|_{pm}, \quad p,m \in \mathbb{N}_0.
    \]
    By Lemma 2.4 from \cite{CGPR2019}, there exists $\delta_0$ such that for every $\delta \geq \delta_0$ we have
    \[
        \sigma_p(\partial_g^\gamma {u}, \delta) \leq 2^p \sigma_0(\partial_g^\gamma {u}, \delta) + \sum_{\ell=1}^{p} 2^{p-\ell} \binom{p}{\ell} (\ell!)^{-\frac{km}{k+m}} \sigma_0(Q^\ell \partial_g^\gamma {u}, \delta).
    \]
    Note that $\sigma_0(\partial_g^\gamma {u}, \delta) = \|\partial_g^\gamma {u}\|_{L^2(\mathbb{R}^n)}$ does not depend on $\delta$. Using the hypothesis \eqref{eq QM M}, we obtain
    \begin{align*}
        \sigma_p(\partial_g^\gamma {u}, \delta) &\leq 2^p C_{\gamma,0} + \sum_{\ell=1}^{p} 2^{p-\ell} \binom{p}{\ell} (\ell!)^{-\frac{km}{k+m}} C_{\ell,\gamma} < \infty,
    \end{align*}
    which concludes the proof.
\end{proof}

\begin{proof}[Proof of Theorem \ref{thm partial schwartz M}]
    Let $u \in \mathcal{S}(\mathfrak{M} \times \mathbb{R}^n)$ and fix $g \in G$ and $\gamma \in \mathbb{N}_0^d$. Then $\partial_g^\gamma u(g,\cdot) \in \mathcal{S}(\mathbb{R}^n)$ and we can write
    \[
        \partial_g^\gamma u(g,x) = \sum_{j \in \mathbb{N}} a_j^\gamma(g) \varphi_j(x),
    \]
    where $a_j^\gamma(g) = \langle \partial_g^\gamma u(g,\cdot), \varphi_j(\cdot) \rangle_{L^2(\mathbb{R}^n)} = \partial_g^\gamma \langle u(g,\cdot), \varphi_j(\cdot) \rangle_{L^2(\mathbb{R}^n)}$. Taking $u_j(g) \doteq a_j^0(g)$, for all $g \in G$, we have
    \[
        \partial_g^\gamma u_j(g) = \partial_g^\gamma a_j^0(g) = a_j^\gamma(g).
    \]
    Now, observe that for any $M \in \mathbb{N}$
    \[
        \|Q^M \partial_g^\gamma u(g,\cdot)\|_{L^2(\mathbb{R}^n)}^2 = \left\| \sum_{j \in \mathbb{N}} \lambda_j^M a_j^\gamma(g) \varphi_j(\cdot) \right\|_{L^2(\mathbb{R}^n)}^2 = \sum_{j \in \mathbb{N}} |\lambda_j|^{2M} |a_j^\gamma(g)|^2.
    \]
    It follows by Proposition \ref{prop 2.3 [AC]22 M} and by the Weyl asymptotics \eqref{weyl asymptotic} that for every $j, M \in \mathbb{N}$ there exists $C_{M,\gamma} > 0$ such that
    \[
        \sup_{g \in G} |a_j^\gamma(g)| \leq C_{M,\gamma} |\lambda_j|^{-M} \leq C_{M,k}' \, j^{-\frac{Mn(m+k)}{mk}}, \quad \forall g \in G.
    \]
    Therefore, given $N \in \mathbb{N}$, taking $M \in \mathbb{N}$ such that $M > \frac{mk N}{n(m+k)}$, we obtain \eqref{eq partial coef uj}.
    
    Conversely, if $u$ is given by
    \[
        u(g,x) = \sum_{j=1}^{\infty} u_j(g) \varphi_j(x),
    \]
    and satisfies \eqref{eq partial coef uj}, then for every $M \in \mathbb{N}$ and $\gamma \in \mathbb{N}_0^d$, by using \eqref{weyl asymptotic}, we have
    \begin{align*}
        \|Q^M \partial_g^\gamma u(g,\cdot)\|_{L^2(\mathbb{R}^n)}^2 &= \left\| \sum_{j \in \mathbb{N}} \lambda_j^M \partial_g^\gamma u_j(g) \varphi_j(\cdot) \right\|_{L^2(\mathbb{R}^n)}^2 \\
        &= \sum_{j \in \mathbb{N}} |\lambda_j|^{2M} |\partial_g^\gamma u_j(g)|^2 \\
        &\leq \sum_{j \in \mathbb{N}} C \, j^{\frac{2Mmk}{n(m+k)}} \, C_{\gamma,N}^2 \, j^{-2N}
    \end{align*}
    for all $g \in G$ and $N \in \mathbb{N}$. Choosing $N > \frac{Mmk}{n(m+k)} + 1$, we obtain that the last series converges. Therefore, by Proposition \ref{prop 2.3 [AC]22 M}, $u \in \mathcal{S}(\mathfrak{M} \times \mathbb{R}^n)$.
\end{proof}

Using the Fourier analysis on $\mathfrak{M}$ described in the previous section, we also obtain a characterization of the Schwartz space on $\mathfrak{M}\times \mathbb{R}^n$ in terms of full Fourier coefficients.

\begin{theorem}\label{thm coef total S M}
    Let $\{\widehat{u}_j(\eta)_{sr}: \ \forall j\in \N, \ \forall [\eta] \in \widehat{G}_0, \  \ 1 \leq r \leq h_\eta,  \ 1 \leq s \leq d_\eta  \}$ be a set of complex numbers and consider the series
    \[ 
        u(g,x) = \sum_{j=1}^{+\infty} \sum_{[\eta] \in \widehat{G}_0} d_\eta \sum_{s=1}^{d_\eta} \sum_{r=1}^{h_\eta} \widehat{u}_j(\eta)_{sr} \, \eta(g)_{rs} \varphi_j(x), \quad (g,x)\in \mathfrak{M} \times \R^n.
    \]
    Then $u\in \mathcal{S}(\mathfrak{M} \times \R^n)$ if and only if for every $N \in \N$ there exists $C_N > 0$ such that
    \begin{equation}\label{coef total S M}
        | \widehat{u}_j(\eta)_{sr} | \leq C_N (j+\langle \eta \rangle )^{-N},\quad \forall j\in \N, \ \forall [\eta] \in \widehat{G}_0, \  \ 1 \leq r \leq h_\eta, \ 1 \leq s \leq d_\eta.
    \end{equation}
\end{theorem}
\begin{proof}
    Assume that \eqref{coef total S M} holds. It follows from Theorem \ref{thm C^infty(G)} that, for each $j\in \N$,
    \[
        u_j = \sum_{[\eta] \in \widehat{G}_0} d_\eta \sum_{s=1}^{d_\eta} \sum_{r=1}^{h_\eta} \widehat{u}_j(\eta)_{sr} \, \eta(g)_{rs}  \in C^\infty(\mathfrak{M}).
    \]    
    Let $ \gamma\in \N_0^d$. 
    From \eqref{eq deta},  \eqref{ex 1.39 [W]}, and \eqref{coef total S M}, for every  $N \in \N$ there exists $C_N > 0$ such that
    \begin{align*}
         |\partial_g^\gamma u_j(g)| &\leq  \sum_{[\eta] \in \widehat{G}_0} d_\eta \sum_{s=1}^{d_\eta} \sum_{r=1}^{h_\eta} |\widehat{u}_j(\eta)_{sr}| \,| \partial_g^\gamma \eta(g)_{rs}|  \\
         &\leq \sum_{[\eta] \in \widehat{G}_0} d_\eta d_\eta h_\eta   C_N (j+\langle \eta \rangle )^{-N}  \, C_\gamma \langle \eta \rangle^M\\
         &\leq \sum_{[\eta] \in \widehat{G}_0} c \langle \eta \rangle^{\frac{d}{2}} d_\eta^2   C_N (j+\langle \eta \rangle )^{-N}  \, C_\gamma \langle \eta \rangle^M,
    \end{align*}
    for some constants $c,C_\gamma,M>0$, where $d$ is the dimension of $G$.
    Then, given $N' \in \N$,  taking $N = N' + 3d+M$, we obtain
     \begin{align*}
         |\partial_g^\gamma u_j(t)|  
         &\leq \sum_{[\eta] \in \widehat{G}_0} c \langle \eta \rangle^{d} d_\eta^2    C_{N'} (j+\langle \eta \rangle )^{-N'}(j+\langle \eta \rangle )^{-3d-M}  \, C_\gamma \langle \eta \rangle^M\\
         &\leq \sum_{[\eta] \in \widehat{G}_0} c \langle \eta \rangle^{d} d_\eta^2  C_{N'} j^{-N'}\langle \eta \rangle ^{-3d-M }  \, C_\gamma \langle \eta \rangle^M\\
         &\leq C_{N', \gamma} j^{-N'} \sum_{[\eta] \in \widehat{G}_0}d_\eta^2 \langle \eta \rangle ^{-2d}.   
    \end{align*}
     Therefore, by \eqref{sum deta < infty}
     and  Theorem \ref{thm partial schwartz M}, $u\in \mathcal{S}( \mathfrak{M}\times \R^n)$.

    On the other hand, if $u \in \mathcal{S}( \mathfrak{M} \times \R^n)$, then,  by Theorem \ref{thm partial schwartz M}, for every $N\in\N$ we have
    \begin{align*}
        \nu_{[\eta]}^N|\widehat{u}_{j}(\eta)_{sr}| &= \left| \widehat{\mathcal{L}_G^N u_j}(\eta)_{sr} \right| = \left| \int_{G} \mathcal{L}_G^N u_j(g) \overline{\eta(g)_{rs}} d\mu_G\right|\\
        &\leq \left( \int_{G} \left| \mathcal{L}_G^N u_j(g) \right|^{2}  d\mu_G \right)^{\frac{1}{2}}  \left( \int_{G} \left| \eta(g)_{rs} \right|^2  d\mu_G \right)^{\frac{1}{2}}\\
        &\leq \dfrac{1}{\sqrt{d_\eta}} \sum_{|\gamma|=2N} \sup_{g\in G}|\partial_g^\gamma u_j(g)|\\
        &\leq  \dfrac{1}{\sqrt{d_\eta}} \sum_{|\gamma|=2N} C_{\gamma,N} j^{-N}.
    \end{align*}
     By \eqref{<eta>^2}, it follows that 
    \begin{align*}
        |\widehat{u}_{j}(\eta)_{sr}| &\leq  
        \dfrac{1}{\sqrt{d_\eta}} \sum_{|\gamma|=2N} C_{\gamma,N} 
         j^{-N}    C^{N} \langle \eta  \rangle^{-2N}\\ &\leq  
        \dfrac{C^{N}}{\sqrt{d_\eta}} \sum_{|\gamma|=2N} C_{\gamma,N}   
         j^{-N}    \langle \eta  \rangle^{-N}  \leq C_N (j+ \langle \eta \rangle)^{-N},
    \end{align*}
    for all $j\in \N$, $[\eta] \in \widehat{G}_0$, $1 \leq r \leq h_\eta$ and $ 1 \leq s \leq d_\eta$, where $C_N \doteq \displaystyle \frac{C^{N}}{\sqrt{d_\eta}} \! \sum_{|\gamma|=2N}\! C_{\gamma,N} 2^{N}$, completing the proof.
\end{proof}

We now present the Fourier analysis of the dual space $\mathcal{S}^\prime(\mathfrak{M} \times \mathbb{R}^n)$ generated by the basis associated with the operator $Q$. The next results are inspired by \cite{AC2022, ACK2025, K2025}.

\begin{theorem}\label{thm parcial1 S' M}
    Let $u\in \mathcal{S}^\prime(\mathfrak{M} \times \mathbb{R}^n) $. For any $j\in\N$, $u_j: C^\infty(\mathfrak{M}) \to \C$ given by
    \[
        \langle u_j,\psi \rangle \doteq \langle u, \psi(g)\varphi_j (x)\rangle, \quad \psi\in C^\infty(\mathfrak{M}),
    \]
    belongs to $\mathcal{D}^\prime(\mathfrak{M} )$, where $g$ represents the variable in $G$ and $x$ the one in $\R^n$, and there exist $C>0$ and $M\in\N$ such that
    \begin{equation}\label{coef parcial S' M}
        |\langle u_j,\psi \rangle|  \leq C \sup_{|\gamma|\leq M}\sup_{g\in G} |\partial_g^\gamma \psi| j^M, \quad \forall j \in \N, \ \psi\in C^\infty(\mathfrak{M}).
    \end{equation}
    Moreover, 
    \begin{equation}
        \langle u, \theta \rangle = \sum_{j=1}^{+\infty} \langle u_j(g) \varphi_j(x), \theta \rangle,
    \end{equation}
    where 
    \[
        \langle u_j(g) \varphi_j(x), \theta \rangle \doteq \left\langle u_j(g) , \int_{\R^n}\theta(g,x)\varphi_j(x) dx \right\rangle,
    \]    
    for all  $\theta\in \mathcal{S}(\mathfrak{M} \times \R^n)$. 
\end{theorem}

\begin{proof}
    Since $u\in \mathcal{S}^\prime(\mathfrak{M} \times \mathbb{R}^n)$, it follows from \eqref{def S' M} that there exist $C,N>0$ such that
    \begin{align*}
        | \langle u_j(g),\psi(g) \rangle|   &= |\langle u, \psi(g)\varphi_j (x)\rangle | \leq C \, p_N(\psi \varphi_j)\\
        &\leq C \sum_{|\gamma|+|\alpha|+|\beta| < N} \sup_{(g,x) \in G \times \mathbb{R}^n} | x^\beta \partial_x^\alpha \partial_g^\gamma  \psi(g)\varphi_j (x)| \\
        &\leq  C \sum_{|\gamma|+|\alpha|+|\beta| < N}  \sup_{g \in G } |  \partial_g^\gamma  \psi(g)| \sup_{x \in \mathbb{R}^n} | x^\beta \partial_x^\alpha  \varphi_j (x)|.
    \end{align*}
    Since  $\varphi_j\in \mathcal{S}(\mathbb{R}^n)$,  the proof of Proposition \ref{prop 2.3 [AC]22 M} shows that the seminorm $\displaystyle \sup_{x \in \mathbb{R}^n} | x^\beta \partial_x^\alpha  \varphi_j (x)|$
    is equivalent to the one in \eqref{eq QM M}. Consequently, there exists $M_0\in \N$ (depending on $\alpha,\beta$) such that
    \begin{align*}
        | \langle u_j(g),\psi(g) \rangle| &\leq  C \sum_{|\gamma|+|\alpha|+|\beta| < N}  \sup_{g \in G } |  \partial_g^\gamma  \psi(g)|  \sum_{M\leq M_0} \| Q^{M} \varphi_j \|_{L^2(\R^n)}\\
        &\leq  C   \sup_{|\gamma| \leq N} \sup_{g \in G } |  \partial_g^\gamma  \psi(g)| \sum_{|\alpha|+|\beta| < N} \sum_{M\leq M_0} |\lambda_j^M|\\
        &\leq C \sup_{|\gamma| \leq N} \sup_{g \in G } |  \partial_g^\gamma  \psi(g)| C_{N} j^{\frac{M_0 mk}{n(m+k)}},
    \end{align*}
   where the last estimate follows from \eqref{weyl asymptotic}, and $C_N> 0$ is a constant depending only on $N$.
    Therefore, { $u_j \in \mathcal{D}^\prime(\mathfrak{M})$ } for all $j\in \N$ and inequality \eqref{coef parcial S' M} holds. Moreover, for $\theta \in \mathcal{S}(\mathfrak{M} \times \mathbb{R}^n)$, we have
    \[
        \theta(g,x) = \sum_{j=1}^{+\infty} \theta_j(g)\varphi_j(x),
    \]    
    where the coefficients are given by $\theta_j(g) = \langle \theta, \varphi_j \rangle = \int_{\R^n} \theta(g,x)\varphi_j(x)dx.$ Then
    \begin{align*}
        \langle u, \theta \rangle &= \sum_{j=1}^{+\infty} \langle u, \theta_j(g)\varphi_j(x) \rangle = \sum_{j=1}^{+\infty} \langle u_j(g), \theta_j(g) \rangle \\
        &=  \sum_{j=1}^{+\infty} \left\langle u_j(g) , \int_{\R^n}\theta(g,x)\varphi_j(x) dx \right\rangle = \sum_{j=1}^{+\infty} \langle u_j(g) \varphi_j(x), \theta \rangle,
    \end{align*}
    completing the proof.
\end{proof}

Conversely, we have the following.

\begin{theorem}\label{thm parcial2 S' M}
    Let $\{u_j\}_{j\in\N} \subset \mathcal{D}'(\mathfrak{M})$ be a sequence for which there exist $C>0$ and $M\in \N$ such that \eqref{coef parcial S' M} holds. Then
    \[
        u(g,x) = \sum_{j=1}^{+\infty} u_j \varphi_j(x).
    \]    
    belongs to $\mathcal{S}'(\mathfrak{M}\times \R^n) $. Moreover,
    \[
        \langle u_j, \psi(g) \rangle = \langle u, \psi(g)\varphi_j (x)\rangle, \quad \psi\in C^\infty(\mathfrak{M}).
    \]    
\end{theorem}

To prove Theorem \ref{thm parcial2 S' M}, we will use the following lemma:

\begin{lemma}\label{Lemma 2.11 [K]25}
    Let $\{u_j\}_{j \in \mathbb{N}} $ be a sequence in $ \mathcal{S}'(\mathfrak{M} \times \mathbb{R}^n)$ such that, for every $\theta \in \mathcal{S}(\mathfrak{M} \times \mathbb{R}^n)$, $\{\langle u_j, \theta \rangle\}_{j \in \mathbb{N}}$ is a Cauchy sequence in $\mathbb{C}$. Then there exists a unique $u \in \mathcal{S}'(\mathfrak{M} \times \mathbb{R}^n)$ such that $u = \displaystyle\lim_{j \to \infty} u_j$ in $\mathcal{S}'(\mathfrak{M} \times \mathbb{R}^n)$.
\end{lemma}
\begin{proof}
    Let $\theta \in \mathcal{S}(\mathfrak{M} \times \mathbb{R}^n)$. Since $\{\langle u_j, \theta \rangle\}_{j \in \mathbb{N}}$ is a Cauchy sequence in $\mathbb{C}$, it converges. We then define 
    \[
        \langle u, \theta \rangle \doteq \lim_{j \to \infty} \langle u_j, \theta \rangle \in \mathbb{C}.
    \]
    Note that the mapping $\theta \mapsto \langle u, \theta \rangle$ is linear, and we claim that it is continuous. Indeed,
    let $\{\theta_\ell\}_{\ell\in\N} \subset \mathcal{S}(\mathfrak{M} \times \mathbb{R}^n) $ be such that $\theta_\ell \to 0$  as $\ell \to \infty$.    
    For each $\psi \in \mathcal{S}(\mathfrak{M} \times \mathbb{R}^n)$, by hypothesis, $\{ \langle u_j, \psi \rangle\}_{j \in \mathbb{N}}$ is a  Cauchy sequence; hence, there exists $C_\psi > 0$ such that $|\langle u_j, \psi \rangle| \leq C_\psi$, for all $j \in \mathbb{N}$. Since $\mathcal{S}(\mathfrak{M} \times \mathbb{R}^n)$ is a Fréchet space and $\mathbb{C}$ is a Banach space, the uniform boundedness principle implies that 
    there exist constants $K > 0$ and $N \in \mathbb{N}$ such that
    \[
        |\langle u_j, \psi \rangle| \leq K 
    \]
    for all $j\in\N$ and  $\psi \in \mathcal{S}(\mathfrak{M} \times \mathbb{R}^n)$ such that $p_N(\psi) \leq 1$.
    Given $\epsilon > 0$, consider $\psi_\ell \doteq \frac{2K}{\epsilon} \theta_\ell$ for all $\ell \in \mathbb{N}$. Therefore, $\psi_\ell \to 0$  in $\mathcal{S}(\mathfrak{M} \times \mathbb{R}^n)$ as $\ell\to +\infty$. Let $\ell_0 \in \mathbb{N}$ be such that $p_N(\psi_\ell) \leq 1$ for all $\ell \geq \ell_0$. It follows that
    \[ 
        |\langle u_j, \theta_\ell \rangle| = \frac{\epsilon}{2K} |\langle u_j, \psi_\ell \rangle| \leq \frac{\epsilon}{2K} \cdot K = \frac{\epsilon}{2}, \quad \forall \ell \geq \ell_0.
    \]
    Moreover, since for each fixed $\ell$, $\langle u, \theta_\ell \rangle = \displaystyle\lim_{j \to \infty} \langle u_j, \theta_\ell \rangle$, there exists $j_\ell \geq \ell$ such that
    \[
        |\langle u, \theta_\ell \rangle - \langle u_{j_\ell}, \theta_\ell \rangle| < \frac{\epsilon}{2}.
    \]    
    Then, for all $\ell \geq \ell_0$,
    \[
        |\langle u, \theta_\ell \rangle| \leq |\langle u, \theta_\ell \rangle - \langle u_{j_\ell}, \theta_\ell \rangle| + |\langle u_{j_\ell}, \theta_\ell \rangle| < \frac{\epsilon}{2} + \frac{\epsilon}{2} = \epsilon,
    \]
    which implies that $\langle u, \theta_\ell \rangle \to 0$ as $\ell \to +\infty$. Therefore, $u \in \mathcal{S}'(\mathfrak{M} \times \mathbb{R}^n)$.
\end{proof}

\begin{proof}[Proof of Theorem \ref{thm parcial2 S' M}]
    Consider for each $\ell \in \N$ the finite sum
    \[
        S_\ell = \sum_{j=1}^\ell u_j \varphi_j(x) \ {\in  \mathcal{S}'(\mathfrak{M}\times \R^n)},
    \]    
    and let $\theta \in \mathcal{S}(\mathfrak{M} \times \mathbb{R}^n)$. For any $k\in \N$ and $\theta \in \mathcal{S}(\mathfrak{M} \times \mathbb{R}^n)$, we have by \eqref{coef parcial S' M} that there exist $C>0$ and $M\in \N $ such that
    \begin{align*}
        \langle S_{\ell+k} - S_\ell , \theta \rangle &\leq \sum_{j=\ell+1}^{\ell+k} |\langle u_j \varphi_j(x), \theta\rangle| = \sum_{j=\ell+1}^{\ell+k} \left|\left\langle u_j  , \int_{\R^n}\theta(g,x)\varphi_j(x) dx \right\rangle \right| \\
        &\leq \sum_{j=\ell+1}^{\ell+k} |\langle u_j , \theta_j(g) \rangle|  \leq \sum_{j=\ell+1}^{\ell+k}  C \sup_{|\gamma|\leq M}\sup_{g\in G} |\partial^\gamma \theta_j| j^M.
    \end{align*}
    Since $\theta \in \mathcal{S}(\mathfrak{M} \times \mathbb{R}^n)$, it follows from Theorem \ref{thm partial schwartz M} that  for every $\gamma\in \N_0^d$ and $N\in \N$ there exists $C_{\gamma,N} > 0$ such that
    \[
        \sup_{g\in G}|\partial^\gamma \theta_j| \leq C_{\gamma,N} j^{-N}.
    \]    
    Then, taking $N= M+2$, we have
    \begin{align*}
         \langle S_{\ell+k} - S_\ell , \theta \rangle &\leq  C \sum_{j=\ell+1}^{\ell+k} \sup_{|\gamma|\leq M} C_{\gamma,M+2} \, j^{-M-2} j^M \leq C_M  \sum_{j=\ell+1}^{\ell+k} j^{-2}.
    \end{align*}
    This implies that $\{\langle S_\ell, \theta \rangle \}_{\ell\in\N} $ is a Cauchy sequence in $\C$, for each $\theta \in \mathcal{S}(\mathfrak{M} \times \mathbb{R}^n)$. Therefore, by Lemma \ref{Lemma 2.11 [K]25}, there exists a unique $u\in \mathcal{S}'(\mathfrak{M} \times \mathbb{R}^n)$ such that
    \[
        u = \lim_{\ell \to +\infty} S_\ell = \sum_{j=1}^{+\infty} u_j \varphi_j,
    \]    
    completing the proof.
\end{proof}

Finally, we characterize the space $\mathcal{S}^\prime(\mathfrak{M} \times \mathbb{R}^n)$ via full Fourier coefficients. 

\begin{theorem}\label{thm coef total S' M}
    Let $\{\widehat{u}_j(\eta)_{sr}: \ \forall j\in \N, \ \forall [\eta] \in \widehat{G}_0, \  \ 1 \leq r \leq h_\eta,  \ 1 \leq s \leq d_\eta  \}$ be a set of complex numbers and consider the series
    \[ 
        u(g,x) = \sum_{j=1}^{+\infty} \sum_{[\eta] \in \widehat{G}_0} d_\eta \sum_{s=1}^{d_\eta} \sum_{r=1}^{h_\eta} \widehat{u}_j(\eta)_{sr} \, \eta(g)_{rs} \varphi_j(x), \quad (g,x)\in G \times \R^n.
    \]    
    Then $u\in \mathcal{S}'(\mathfrak{M} \times \R^n)$ if and only there exist $C>0$ and $M \in \N$ such that
    \begin{equation}\label{coef total S' M}
        | \widehat{u}_j(\eta)_{sr} | \leq C (j+\langle \eta \rangle )^{M},\quad \forall j\in \N, \ \forall [\eta] \in \widehat{G}_0, \  \ 1 \leq r \leq h_\eta, \ 1 \leq s \leq d_\eta.
    \end{equation}
\end{theorem}
\begin{proof}
    Assume that \eqref{coef total S' M} holds. It follows from Theorem \ref{thm D'(G)}     that, for each $j\in \N$,
    \[
        u_j = \sum_{[\eta] \in \widehat{G}_0} d_\eta \sum_{s=1}^{d_\eta} \sum_{r=1}^{h_\eta} \widehat{u}_j(\eta)_{sr} \, \eta(g)_{rs}
    \]
    belongs to $\mathcal{D}'(\mathfrak{M})$.  Consider for each $\ell \in \N$ the sum
    \[
        S_\ell = \sum_{j=1}^\ell u_j \varphi_j(x) \ {\in  \mathcal{S}'(\mathfrak{M} \times \R^n).}
    \]
    For any $k\in \N$ and $\theta \in \mathcal{S}(\mathfrak{M} \times \R^n) $, we have
    \begin{align*}
    \langle S_{\ell+k} - S_\ell , \theta \rangle &= \sum_{j=\ell+1}^{\ell+k} \langle u_j \varphi_j(x), \theta\rangle = \sum_{j=\ell+1}^{\ell+k} \langle u_j , \theta_j \rangle\\
    &= \sum_{j=\ell+1}^{\ell+k}  \sum_{[\eta] \in \widehat{G}_0} d_\eta \sum_{s=1}^{d_\eta} \sum_{r=1}^{h_\eta} \widehat{u}_j(\eta)_{sr} \langle   \eta(g)_{rs} , \theta_j \rangle\\
    &=  \sum_{j=\ell+1}^{\ell+k}  \sum_{[\eta] \in \widehat{G}_0} d_\eta \sum_{s=1}^{d_\eta} \sum_{r=1}^{h_\eta} \widehat{u}_j(\eta)_{sr} \int_G    \eta(g)_{rs} \, \theta_j(g) d\mu_G\\
    &= \sum_{j=\ell+1}^{\ell+k}  \sum_{[\eta] \in \widehat{G}_0} d_\eta \sum_{s=1}^{d_\eta} \sum_{r=1}^{h_\eta} \widehat{u}_j(\eta)_{sr} \int_G    \theta_j(g)\, \overline{\overline{\eta}(g)_{rs}} d\mu_G\\
    &= \sum_{j=\ell+1}^{\ell+k}  \sum_{[\eta] \in \widehat{G}_0} d_\eta \sum_{s=1}^{d_\eta} \sum_{r=1}^{h_\eta} \widehat{u}_j(\eta)_{sr} \widehat{\theta}_j(\overline{\eta})_{sr}
\end{align*}
    Since $\theta \in \mathcal{S}(\mathfrak{M} \times \R^n)$, by Theorem \ref{thm coef total S M} for every $N\in \N$ there exist $C_N > 0$ such that 
    \[
        | \widehat{\theta}_j(\overline{\eta})_{sr}|  \leq C_N (j+\langle \eta \rangle )^{-N}.
    \]    
    Taking $N= 3d + 2M$, by hypothesis \eqref{coef total S' M} we obtain 
    \begin{align*}
        |\langle S_{\ell+k} - S_\ell , \theta \rangle| &\leq    \sum_{j=\ell+1}^{\ell+k}  \sum_{[\eta] \in \widehat{G}_0} d_\eta \, d_\eta^ 2  \, C(j+\langle \eta \rangle )^{M}\, C_{N}(j+\langle \eta \rangle )^{-2M-3d}\\
        &\leq C_M   \sum_{j=\ell+1}^{\ell+k}  \sum_{[\eta] \in \widehat{G}_0}  c \langle \eta \rangle^{d} \, d_\eta^2 \, j^{-M} \langle \eta \rangle^{-3d}\\
        &\leq C_M' \, j^{-M}\sum_{j=\ell+1}^{\ell+k}  \sum_{[\eta] \in \widehat{G}_0} d_\eta^2 \langle \eta \rangle^{-2d},
    \end{align*}
    using \eqref{eq deta}, and the last series converges by \eqref{sum deta < infty}.
    Therefore, $\{\langle S_\ell, \theta \rangle \}_{\ell\in\N} $ is a Cauchy sequence in $\C$, for all $\theta \in \mathcal{S}(\mathfrak{M} \times \mathbb{R}^n)$ and, by Lemma \ref{Lemma 2.11 [K]25}, there exists $u\in \mathcal{S}'(\mathfrak{M} \times \mathbb{R}^n)$ such that $S_\ell \to u $ in $\mathcal{S}'(\mathfrak{M} \times \mathbb{R}^n)$ as $\ell \to +\infty$.

    On the other hand, let $u\in \mathcal{S}'(\mathfrak{M}\times \R^n)$. By Theorems \ref{thm parcial1 S' M} and \ref{thm parcial2 S' M}, we can write 
    \[
        u(t,x) = \sum_{j=1}^{+\infty} u_j \varphi_j(x),
    \]
    where $u_j \in \mathcal{D}'(\mathfrak{M})$ is given by $\langle u_j, \psi(g) \rangle = \langle u, \psi(g)\varphi_j (x)\rangle$, for all $\psi\in C^\infty(\mathfrak{M})$. Therefore,
    by Theorem \ref{thm D'(G)},
    \[
        u(t,x) = \sum_{j=1}^{+\infty} \sum_{[\eta] \in \widehat{G}_0} d_\eta \sum_{s=1}^{d_\eta} \sum_{r=1}^{h_\eta} \widehat{u}_j(\eta)_{sr} \, \eta(g)_{rs}  \varphi_j(x),
    \]
     with $ \widehat{u}_j(\eta)_{sr} = \langle u_j , \overline{\eta_{rs}} \rangle_G $.
     By \eqref{ex 1.39 [W]} and  \eqref{coef parcial S' M}, there exists $C>0$ and $M\in \N$ such that
     \begin{align*}
         |\widehat{u}_j(\eta)_{sr}| &=  |\langle u_j , \overline{\eta_{rs}} \rangle|  \leq  C \sup_{|\gamma|\leq M}\sup_{g\in G} |\partial_g^\gamma \overline{\eta_{rs}}(g)   | j^M\\
         &\leq C_M \langle \eta \rangle^M j^M \leq C(j+\langle \eta \rangle)^{2M},
     \end{align*}
     completing the proof.
\end{proof}

\section{\texorpdfstring{Schwartz regularity and solvability of a class of operators on $\mathfrak{M} \times \R^n$}{Schwartz regularity and solvability of a class of operators on M x Rn}}

Let $Q$ the partial differential operator given by \eqref{def Q} satisfying \eqref{anistropic eliptic}. Let $T$ a left-invariant and symmetric operator on $G$, and we consider $T f \doteq T \tilde{f}$ for all $f\in C^\infty(\mathfrak{M})$. For each  $ [\eta] \in \widehat{G}_0$, we can choose a representative  $\eta$  such that  $ \sigma_T $ is a diagonal matrix  
\[
    \sigma_T(\eta)_{rs} = \mu_{T,r}(\eta) \delta_{rs}, \quad 1 \leq r,  s \leq d_\eta,
\]
where $  \mu_{T,r}(\eta) \in \mathbb{R}$, and we can write
\[
    T\eta(g)_{rs} =  \mu_{T,s}(\eta)\eta(g)_{rs}, \quad 1 \leq r \leq h_\eta, \ 1\leq s \leq d_\eta, \ \forall g\in \mathfrak{M}.
\]

Consider on $\mathfrak{M} \times \R^n$ the operator
\begin{equation}\label{op P=T+cQ+q}
     \mathsf{P} = T + c Q + q,
\end{equation}
where $c, \ q \in \C$, $c\neq 0$. 
Our goal is to  study the Schwartz regularity of $\mathsf{P}$ in the following sense: 

\begin{definition}
    We say that the operator $\mathsf{P}$ is $\mathcal{S}$-globally hypoelliptic on $ \mathfrak{M}\times \R^n$ if 
    \[
        u\in \mathcal{S}^\prime(\mathfrak{M} \times \R^n)  \text{ and }  \mathsf{P} u \in   \mathcal{S}(\mathfrak{M} \times \R^n)  \implies u \in \mathcal{S}(\mathfrak{M} \times \R^n).
    \]    
\end{definition}

Let $u\in \mathcal{S}^\prime(\mathfrak{M} \times \R^n)$ be a solution for $\mathsf{P} u=f \in \mathcal{S}'(\mathfrak{M} \times \R^n)$. Expanding $u$ and $f$ in their total Fourier series, we have by Theorems \ref{thm coef total S M} and \ref{thm coef total S' M} that
\[
    u(g,x) = \sum_{j=1}^{+\infty} \sum_{[\eta] \in \widehat{G}_0} d_\eta \sum_{s=1}^{d_\eta} \sum_{r=1}^{h_\eta} \widehat{u}_j(\eta)_{sr} \, \eta(g)_{rs} \varphi_j(x), 
\]
and 
\[
    f(g,x) = \sum_{j=1}^{+\infty} \sum_{[\eta] \in \widehat{G}_0} d_\eta \sum_{s=1}^{d_\eta} \sum_{r=1}^{h_\eta} \widehat{f}_j(\eta)_{sr} \, \eta(g)_{rs} \varphi_j(x),
\]
for all $ (g,x)\in \mathfrak{M} \times \R^n$.
This implies that  
\begin{equation}\label{eq coef}
     {(\mu_{T,s}(\eta) + c \lambda_j + q) \widehat{u}_j(\eta)_{sr} =  \widehat{f}_j(\eta)_{sr},}
\end{equation}
for all $ j \in \N $, $[\eta]\in \widehat{G}_0 $, and $ 1 \leq r \leq h_\eta, \ 1 \leq s \leq d_\eta$.
Define the set
\[
    \mathcal{Z}= \{ (j,[\eta])\in \N\times \widehat{G}_0: \  \mu_{T,s}(\eta) + c \lambda_j + q = 0, \quad \text{for some } 1 \leq s \leq d_\eta \}.
\]
\begin{prop}\label{prop Z infinito M}
    If the set $\mathcal{Z}$ has  infinitely many elements, then the operator $\mathsf{P}$ in \eqref{op P=T+cQ+q} is not  $\mathcal{S}$-globally hypoelliptic on $\mathfrak{M} \times \R^n$. 
\end{prop}
\begin{proof}
    Define 
    \[  
         \widehat{u}_j(\eta)_{sr} = \begin{cases}
        1, &\text{if }\  \mu_{T,s}(\eta) + c \lambda_j + q = 0, \\
            0, &\text{otherwise},
        \end{cases}
    \]
    for all $(j,[\eta])\in \N\times \widehat{G}_0$,  and $ 1 \leq r \leq h_\eta, \ 1 \leq s \leq d_\eta$. Then, by Theorem \ref{thm coef total S M} and Theorem \ref{thm coef total S' M}, since $\mathcal{Z}$ has  infinitely many elements,
    \[
        u(g,x) = \sum_{j=1}^{+\infty} \sum_{[\eta] \in \widehat{G}_0} d_\eta \sum_{s=1}^{d_\eta} \sum_{r=1}^{h_\eta} \widehat{u}_j(\eta)_{sr} \, \eta(g)_{rs} \varphi_j(x) \in \mathcal{S}^\prime(\mathfrak{M} \times \R^n) \setminus \mathcal{S}(\mathfrak{M} \times \R^n).
    \]
    Since $ \mathsf{P} u = 0$, we conclude that $ \mathsf{P}$ is not $\mathcal{S}$-globally hypoelliptic on $\mathfrak{M}\times \R^n$.  
\end{proof}

Now let us establish a condition useful to characterize the regularity of $\mathsf{P}$.

\begin{definition}
    We say that $\mathsf{P}$ satisfies the algebraic condition \eqref{AC} if there exist $C>0$ and $M\in \N$ such that
    \begin{equation}\label{AC}\tag{AC}
        | \mu_{T,s}(\eta) + c \lambda_j + q| \geq C (j+\langle \eta \rangle)^{-M},
    \end{equation}
for all $ (j,[\eta])\in \N\times \widehat{G}_0$ and $ 1 \leq s \leq d_\eta$ such that $| \mu_{T,s}(\eta) + c \lambda_j + q|\neq 0$.
\end{definition}

\begin{theorem}\label{thm hip cte}
    The operator $\mathsf{P}$ in \eqref{op P=T+cQ+q} is $\mathcal{S}$-globally hypoelliptic on $\mathfrak{M} \times \R^n$ if and only if the set $\mathcal{Z}$ is finite and $\mathsf{P}$ satisfies \eqref{AC}.
\end{theorem}
\begin{proof}
    By Proposition \ref{prop Z infinito M}, if $\mathcal{Z}$ has  infinitely many elements, then $\mathsf{P}$ is not $\mathcal{S}$-globally hypoelliptic on $\mathfrak{M} \times \R^n$. Assume now that $\mathcal{Z}$ is finite and $\mathsf{P}$ does not satisfy \eqref{AC}. Then, for each $ \ell \in \N$ there exists  $(j_\ell, [\eta_\ell ]) \in \mathbb{N}\times\widehat{G}_0$ such that
    \begin{equation}\label{eq nao AC}
        0<|\mu_{T,s}(\eta_\ell) + c \lambda_{j_\ell} + q| < \frac{1}{\ell}\, (j_\ell + \langle \eta_\ell \rangle)^{-\ell}, \quad \forall \ell \in \mathbb{N}.
    \end{equation}
    Taking a subsequence if necessary, we can choose the pairs $(j_\ell,[\eta_\ell])$ such that $ (j_\ell,[\eta_\ell]) \neq (j_k,[\eta_k])$ for every $\ell \neq k$ in $\N$.
    Indeed, suppose that the pairs $(j_\ell, [\eta_\ell])$ take only finitely many distinct values. Then, there exist $\ell_0 \in \N$ and an infinite set $\Lambda \subset \N$ such that $(j_\ell, [\eta_\ell]) = (j_{\ell_0}, [\eta_{\ell_0}])$, for all $\ell \in \Lambda$. This implies that $|\mu_{T,s}(\eta_\ell) + c \lambda_{j_\ell} + q| = |\mu_{T,s}(\eta_{\ell_0}) + c \lambda_{j_{\ell_0}} + q|, $ for all $\ell \in \Lambda$, and hence it does not converge to zero as $\ell \to \infty$, $\ell \in \Lambda$,  which contradicts inequality \eqref{eq  nao AC}.
    
    Define 
    \[
        f(g,x) \doteq  \sum_{\ell\in \N} d_{\eta_\ell} \sum_{s=1}^{d_{\eta_\ell}} \sum_{r=1}^{h_{\eta_\ell}} (\mu_{T,s}(\eta_\ell) + c \lambda_{j_\ell} + q)  \eta_\ell(g)_{rs} \varphi_{j_\ell}(x).
    \]
    Then, by \eqref{eq nao AC} and Theorem \ref{thm coef total S M}, we have $f\in \mathcal{S}(\mathfrak{M} \times \R^n)$.
    Now, define
    \[
        u(g,x) \doteq  \sum_{\ell\in \N} d_{\eta_\ell} \sum_{s=1}^{d_{\eta_\ell}} \sum_{r=1}^{h_{\eta_\ell}}  \eta_\ell(g)_{rs} \varphi_{j_\ell}(x).
    \]
    By Theorem \ref{thm coef total S' M}, we have $u \in \mathcal{S}'(\mathfrak{M} \times \R^n) \setminus \mathcal{S}(\mathfrak{M} \times \R^n)$.
    Moreover, by construction,
    $\mathsf{P}u=f$. This concludes the proof that $\mathsf{P}$ is not $\mathcal{S}$-globally hypoelliptic  on $\mathfrak{M} \times \R^n$.
    
    Conversely, assume that $\mathcal{Z}$ is finite and $\mathsf{P}$ satisfies \eqref{AC}. 
    If $u\in \mathcal{S}^\prime(\mathfrak{M} \times \R^n)$ is a solution for $\mathsf{P}u=f \in \mathcal{S}(\mathfrak{M} \times \R^n)$, it must satisfy
    \[
        \widehat{u}_j(\eta)_{sr} = \dfrac{\widehat{f}_j(\eta)_{sr}}{ \mu_{T,s}(\eta) + c \lambda_j + q},
    \]
    for all $ (j,[\eta]) \in  \N\times \widehat{G}_0 \setminus \mathcal{Z}$, and $ 1 \leq r \leq h_\eta, \ 1 \leq s \leq d_\eta$.
    By Theorem \ref{thm coef total S M} and condition \eqref{AC}, there exist $C>0$,  $M\in \N$ and for every $N\in \N$, there exist $C_N$ such that
    \begin{align*}
        |\widehat{u}_j(\eta)_{sr}| &= \dfrac{|\widehat{f}_j(\eta)_{sr}|}{| \mu_{T,s}(\eta) + c \lambda_j + q| } \leq C_NC^{-1} (j+\langle \eta \rangle)^{-N + M}, 
    \end{align*}
    for all $ (j,[\eta]) \in  \N\times \widehat{G}_0 \setminus \mathcal{Z}$ and $ 1 \leq r \leq h_\eta, \ 1 \leq s \leq d_\eta$. Given $N'\in \N$, let $N= N' + M$. Since $\mathcal{Z}$ is finite, we may increase the constant if necessary to obtain $C_{N'}> 0$ such that
    \[
        |\widehat{u}_j(\eta)_{sr}| \leq C_{N'} (j+\langle \eta \rangle)^{-N '}, \quad \forall j\in \N, \ \forall [\eta] \in \widehat{G}_0, \  \ 1 \leq r \leq h_\eta, \ 1 \leq s \leq d_\eta.
    \]
    By Theorem \ref{thm coef total S M}, $ u\in \mathcal{S}(\mathfrak{M} \times \R^n)$, which means that $\mathsf{P}$ is $\mathcal{S}$-globally hypoelliptic.
\end{proof}

Now we proceed with the study of the solvability of the operator $\mathsf{P}$. Let $u\in \mathcal{S}'(\mathfrak{M}\times \R^n)$ be a solution of the equation $\mathsf{P}u=f$, where $f\in \mathcal{S}'(\mathfrak{M}\times \R^n)$. By equation \eqref{eq coef},
if 
\[
    \mu_{T,s}(\eta) + c \lambda_j + q = 0,
\]
for some $ j \in \N $, $[\eta]\in \widehat{G}_0 $, and $ 1 \leq s \leq d_\eta$, then $ \widehat{f}_j(\eta)_{sr} = 0$ for every $ 1 \leq r \leq h_\eta$. In view of this necessary condition, we introduce the following compatibility sets:
\[
    \mathcal{F} = \{ f \in \mathcal{S}(\mathfrak{M}\times \R^n): \  \widehat{f}_j(\eta)_{sr} = 0, \text{ whenever } \mu_{T,s}(\eta) + c \lambda_j + q =0\},
\] 
and
\[
    \mathcal{F}^{\prime} = \{ f \in \mathcal{S}'(\mathfrak{M}\times \R^n): \  \widehat{f}_j(\eta)_{sr} = 0, \text{ whenever } \mu_{T,s}(\eta) + c \lambda_j + q =0\}.
\]
Accordingly, we consider the following notions of solvability.
\begin{definition}
    The operator $\mathsf{P}$ is said to be:
    \begin{enumerate}[label=\arabic*.]
        \item 
        $\mathcal{S}$-globally solvable on  $\mathfrak{M}\times \R^n$ if $ \mathsf{P}(\mathcal{S}(\mathfrak{M}\times \R^n) ) = \mathcal{F}$. 
        \item $\mathcal{S}'$-globally solvable on  $\mathfrak{M}\times \R^n$ if $ \mathsf{P}(\mathcal{S}'(\mathfrak{M}\times \R^n) ) = \mathcal{F}'$. 
    \end{enumerate}       
\end{definition}

The next result shows that solvability can also be characterized by condition \eqref{AC}.

\begin{theorem}\label{thm solv cte}
    The operator $\mathsf{P}$ in \eqref{op P=T+cQ+q} is $\mathcal{S}$-globally solvable on  $\mathfrak{M}\times \R^n$ if and only if $\mathsf{P}$ satisfies the condition \eqref{AC}.
\end{theorem}
\begin{proof}
    First, suppose that $\mathsf{P}$ satisfies the condition \eqref{AC} and let $f\in \mathcal{F}$. Define 
    \[  
         \widehat{u}_j(\eta)_{sr} = \begin{cases}
        0, &\text{if }\  \mu_{T,s}(\eta) + c \lambda_j + q = 0, \\
           (\mu_{T,s}(\eta) + c \lambda_j + q)^{-1} \widehat{f}_j(\eta)_{sr}, &\text{otherwise},
        \end{cases}
    \]
    for all $(j,[\eta])\in \N\times \widehat{G}_0$,  and $ 1 \leq r \leq h_\eta, \ 1 \leq s \leq d_\eta$. 
    Since $f\in \mathcal{S}(\mathfrak{M}\times \R^n)$, by Theorem \ref{thm coef total S M}, for every $N\in \N$, there exists ${C}_N>0$ such that
    \[
        | \widehat{f}_j(\eta)_{sr} | \leq C_N (j+\langle \eta \rangle )^{-N},\quad \forall j\in \N, \ \forall [\eta] \in \widehat{G}_0, \  \ 1 \leq r \leq h_\eta, \ 1 \leq s \leq d_\eta. 
    \] 
    By condition \eqref{AC}, there exist $C>0$ and $M\in\N$ such that 
    \begin{align*}
        |\widehat{u}_j(\eta)_{sr}| &\leq |\mu_{T,s}(\eta) + c \lambda_j + q |^{-1} |\widehat{f}_j(\eta)_{sr}| \leq C^{-1} C_N (j+\langle \eta \rangle )^{M-N},
    \end{align*}
    for all $(j,[\eta])\in \N\times \widehat{G}_0$,  and $ 1 \leq r \leq h_\eta, \ 1 \leq s \leq d_\eta$. Therefore, $u\in \mathcal{S}(\mathfrak{M}\times \R^n)$, which implies that  $\mathsf{P}$ is $\mathcal{S}$-globally solvable on  $\mathfrak{M}\times \R^n$.

    Conversely, suppose that $\mathsf{P}$ does not satisfy \eqref{AC}. Hence, there exists a sequence $\{(j_\ell, [\eta_\ell ])\}_{\ell\in \N} \subset \mathbb{N}\times\widehat{G}_0$ of pairwise distinct elements such that
    \begin{equation}\label{eq nao AC 2}
        0<|\mu_{T,s}(\eta_\ell) + c \lambda_{j_\ell} + q| < \dfrac{1}{\ell} \, (j_\ell + \langle \eta_\ell \rangle)^{-\ell}, \quad \forall \ell \in \mathbb{N}.
    \end{equation}
    Consider
    \[
        f(g,x) =  \sum_{\ell\in \N} d_{\eta_\ell} \sum_{s=1}^{d_{\eta_\ell}} \sum_{r=1}^{h_{\eta_\ell}}  ( j_\ell + \langle \eta_\ell \rangle)^{-\frac{\ell}{2}}\eta_\ell(g)_{rs} \varphi_{j_\ell}(x) \in \mathcal{F}.
    \]
    If there exists a solution $u$ for  $\mathsf{P}u=f$, by equation \eqref{eq coef}, the coefficients of $u$ must satisfy
    \begin{align*}
        |\widehat{u}_{j_\ell}(\eta_\ell)_{sr} | = |\mu_{T,s}(\eta_\ell) + c \lambda_{j_\ell} + q|^{-1} (j_\ell + \langle \eta_\ell \rangle)^{-\frac{\ell}{2}} > (j_\ell + \langle \eta_\ell \rangle)^{\frac{\ell}{2}},
    \end{align*}
    for all $\ell \in \N$, $ 1 \leq r \leq h_{\eta_\ell}$, and $1 \leq s \leq d_{\eta_\ell}$. Therefore, by Theorem \ref{thm coef total S' M}, $u\notin \mathcal{S}'(\mathfrak{M}\times \R^n)$ and, therefore, $\mathsf{P} $  is not     $\mathcal{S}$-globally solvable on  $\mathfrak{M}\times \R^n$.
\end{proof}

Observe that the same arguments used in the proof above show that the condition
\eqref{AC} is also necessary and sufficient for $\mathcal{S}^\prime$-global solvability on $\mathfrak{M}\times\mathbb{R}^n$. Indeed, in the first part of the proof, if $f\in\mathcal{F}'$, it is enough to apply Theorem~\ref{thm coef total S' M} to conclude that $u\in\mathcal{S}'(\mathfrak{M}\times\mathbb{R}^n)$. Hence, both notions of global solvability are equivalent:

\begin{corollary}
    The operator $\mathsf{P}$ is 
    $\mathcal{S}$-globally solvable on $\mathfrak{M} \times \R^n$ if and only if it is $\mathcal{S}'$-globally solvable on $\mathfrak{M} \times \R^n$.  
\end{corollary}

Also, as a consequence of the proof of Theorem \ref{thm solv cte} we obtain:
\begin{theorem}
    The operator $\mathsf{P}$ satisfies the condition \eqref{AC} if and only if for all $f\in \mathcal{F}$ there exists $u\in \mathcal{S}^\prime(\T^1\times \R^n)$ such that $\mathsf{P}u=f$.
\end{theorem}

Moreover, by combining Theorems \ref{thm hip cte} and \ref{thm solv cte}, we obtain the following result.

\begin{corollary}
     If the operator $\mathsf{P}$ is $\mathcal{S}$-globally hypoelliptic on $\mathfrak{M} \times \R^n$, then it is also $\mathcal{S}$-globally solvable on  $\mathfrak{M}\times \R^n$. 
\end{corollary}

Solvability is also closely related in the literature to the notion of almost global hypoellipticity (see, for instance, \cite{CDRS2026}). In our setting, we therefore consider the following definition.

\begin{definition}
    We say that the operator $\mathsf{P}$ is almost $\mathcal{S}$-globally hypoelliptic on $\mathfrak{M}\times \R^n$ if, for every $u\in\mathcal{S}^\prime(\mathfrak{M}\times\mathbb{R}^n)$ such that $\mathsf{P}u\in \mathcal{S}(\mathfrak{M}\times \R^n)$, there exists $v\in \mathcal{S}(\mathfrak{M}\times \R^n)$ such that $ \mathsf{P}v=\mathsf{P}u$.
\end{definition}

\begin{theorem}
    The operator $\mathsf{P}$ is     $\mathcal{S}$-globally solvable on  $\mathfrak{M}\times \R^n$ if and only if $\mathsf{P}$ is  almost $\mathcal{S}$-globally hypoelliptic on $\mathfrak{M}\times \R^n$. 
\end{theorem}
\begin{proof}
    Assume that $\mathsf{P}$ is $\mathcal{S}$-globally solvable. Let $u\in\mathcal{S}^\prime(\mathfrak{M}\times\mathbb{R}^n)$ be such that $\mathsf{P}u\in \mathcal{S}(\mathfrak{M}\times \R^n)$. It follows directly from the solvability of $\mathsf{P}$ that there exists $v\in \mathcal{S}(\mathfrak{M}\times \R^n)$ such that $\mathsf{P} v= \mathsf{P} u $.

    Conversely, suppose that $\mathsf{P}$ is not $\mathcal{S}$-globally solvable. By, Theorem \ref{thm solv cte}, $\mathsf{P}$ does not satisfy the condition \eqref{AC}. Hence, there exists a sequence $\{(j_\ell, [\eta_\ell ])\}_{\ell\in \N} \subset \mathbb{N}\times\widehat{G}_0$ of pairwise distinct elements such that
    \begin{equation}\label{2 eq nao AC 2}
        0<|\mu_{T,s}(\eta_\ell) + c \lambda_{j_\ell} + q| < \dfrac{1}{\ell} \, (j_\ell + \langle \eta_\ell \rangle)^{-\ell}, \quad \forall \ell \in \mathbb{N}.
    \end{equation} 
    Consider the function 
    \[
        f(g,x) \doteq  \sum_{\ell\in \N} d_{\eta_\ell} \sum_{s=1}^{d_{\eta_\ell}} \sum_{r=1}^{h_{\eta_\ell}} (\mu_{T,s}(\eta_\ell) + c \lambda_{j_\ell} + q)  \eta_\ell(g)_{rs} \varphi_{j_\ell}(x).
    \]
    By \eqref{2 eq nao AC 2} and Theorem \ref{thm coef total S M}, we have $f\in \mathcal{F}$.
    Observe that $u \in \mathcal{S}^\prime(\mathfrak{M}\times \R^n)$ given by the coefficients
    \[  
        \widehat{u}_j(\eta)_{sr} = \begin{cases}
        1 , &\text{if } j=j_\ell \text{ and } \eta=\eta_\ell, \\
         0, &\text{otherwise},
        \end{cases}
    \]
    is a solution for $\mathsf{P}u=f$. However, by equation \eqref{eq coef}, the coefficients of any solution  $v$ for the equation $\mathsf{P}v=f$, must satisfy
    \begin{align}\label{eq coef ujl}
        |\widehat{v}_{j_\ell}(\eta_\ell)_{sr} | = |\mu_{T,s}(\eta_\ell) + c \lambda_{j_\ell} + q|^{-1} |\mu_{T,s}(\eta_\ell) + c \lambda_{j_\ell} + q| =1,
    \end{align}
    for all $\ell \in \N$, $ 1 \leq r \leq h_{\eta_\ell}$, and $1 \leq s \leq d_{\eta_\ell}$. It implies that $v\notin \mathcal{S}(\mathfrak{M}\times \R^n)$. Therefore, $\mathsf{P}$ is not almost $\mathcal{S}$-globally hypoelliptic on $\mathfrak{M}\times \R^n$.
\end{proof}

\section{Acknowledgements}
I. V. Coelho da Silva was financed by the Coordenação de Aperfeiçoamento de Pessoal de Nível Superior – Brazil (CAPES) – Finance Code 001. P. L. Dattori da Silva was supported in part by São Paulo Research Foundation (FAPESP), Brazil, Process Numbers $\#$2024/08416-6 and $\#$2024/12753-8, and National Council for Scientific and Technological Development, CNPq, Brazil (grant $\#$312210/2026-4). M. Ruzhansky was supported by the Methusalem programme of the Ghent University Special Research Fund (grant 01M01021) and FWO-FAPESP Scientific Cooperation Research (grant G.0AOZ.25N).

\bibliographystyle{elsarticle-num} 
  \bibliography{references.bib}
\end{document}